\documentclass[12pt]{article}
\usepackage{amsmath,amsthm,amssymb,epsfig,graphics,subfigure,color}
\usepackage{tikz}

\input amssym.def
\input amssym.tex

\begin{document}

\centerline{ {\Large \bf Fifty two years ago in Jerusalem}}
\vspace{5mm}

\centerline {{ \bf Lech Maligranda (Lule\aa)}}
\vspace{5mm}

\noindent
On 5-12 July 1960 in Jerusalem was held {\it International Symposium on Linear Spaces} 
at the Hebrew University of Jerusalem under auspicies of the International Mathematical 
Union. The nice photo of 47 speakers and participants was done there. Larger audience 
of mathematicians should see it since we can find here several famous mathematicians.
\vspace{1mm}

\begin{figure}[th!]

\centerline{\includegraphics[width=16cm,height=12cm]{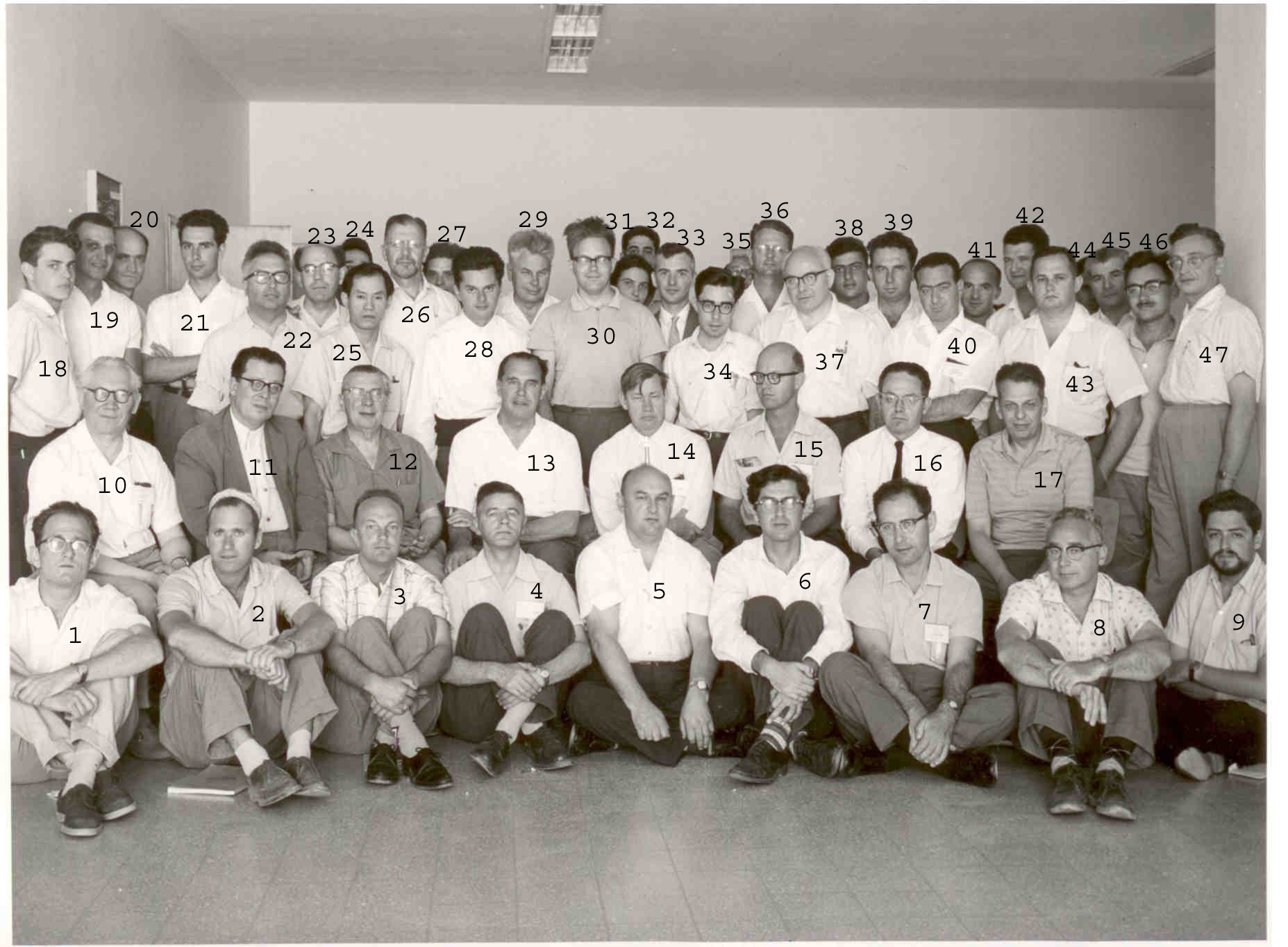}}

\begin{minipage}{160mm}

\end{minipage}
\end{figure}
\vspace{-8mm}

{\footnotesize {\bf Jerusalem, July 1960.} 1. John {\bf Wermer} (1927), 2. Leon {\bf Ehrenpreis} (1930--2010), 
3. Wilhelmus Anthonius Josephus {\bf Luxemburg} (1929), 4. Adriaan Cornelis {\bf Zaanen} (1913--2003), 
5. Aryeh {\bf Dvoretzky} (1916--2008), 6. Michael Bahir {\bf Maschler} (1927--2008), 7. Israel {\bf Halperin} 
(1911--2007), 8. Ralph Saul {\bf Phillips} (1913--1998), 9. Louis {\bf Nirenberg} (1925), 10. W{\l}adys{\l}aw 
{\bf Orlicz} (1903--1990), 11. Gottfried {\bf K\"othe} (1905--1989), 
12. Einar {\bf Hille} (1894--1980), 13. Jean Alexandre Eug\`ene {\bf Dieudonn\'e} (1906--1992), 14. Marshall 
Harvey {\bf Stone} (1903--1989), 15. Richard {\bf Arens} (1919--2000), 16. George Whitelaw {\bf Mackey} 
(1916--2006), 17. Lipman {\bf Bers} (1914--1993), 
18. Daniel {\bf Sternheimer} (1938), 19. unknown, 20. Irving {\bf Kaplansky} 
(1917--2006), 21. Jean-Pierre {\bf Kahane} (1926), 22. Jos\'e Luis {\bf Massera} (1915--2002), 23. Abraham 
{\bf Robinson} (1918--1974), 24. Eliahu {\bf Shamir} (1934), 25. Ichiro {\bf Amemiya} (1923--1995), 26. Angus 
Ellis {\bf Taylor} (1911--1999), 27. Dan {\bf Amir} (1933), 28. Emilio {\bf Gagliardo} (1930--2008), 29. Robert 
Elston {\bf Fullerton} (1906-1963), 30. Victor {\bf Klee} (1925--2007), 31. Tamar {\bf Berger-Burak} (1938), 
32. unknown, 33. Henry {\bf Helson} (1927--2010), 34. Harry {\bf Kesten} (1931), 35. Haim {\bf Amsterdamer}, 
36. Frank Featherston {\bf Bonsall} (1920-2011), 37. Nachman {\bf Aronszajn} (1907--1980), 38. Joram 
{\bf Lindenstrauss} (1936-2012), 39. Shmuel {\bf Agmon} (1922), 40. Guido {\bf Stampacchia} (1922--1978), 
41. Amram {\bf Meir} (1929), 42. Peter David {\bf Lax} (1926), 43. Leopoldo {\bf Nachbin} (1922--1993), 
44. unknown, 45. Meir {\bf Reichaw} [Marian {\bf Reichbach}] (1923-2000), 46. Edmond Ernest {\bf Granirer} (1935), 
47. Paul {\bf Katz} (1924-2005).} \\

\vspace{-5mm}
{\small Three mathematicians are not recognized on the photo, namely, number 19, 32 and 44 . I will be grateful 
if anyone can recognize them and inform me about this. I would like to mention that two more mathematicians were 
at the symposium but I don't see them on the photo: Gaetano {\bf Fichera} (1922--1996) and Jan {\bf Mikusinski} 
(1913--1987).}
\vspace{1mm}

Lectures were given on Tuesday, July 5 by Stone, Dieudonn\'e, Aronszajn, Klee, Dvoretzky, Bers; on Wednesday, 
July 6 by Ehrenpreis, Fichera, Hille, Taylor, Massera; on Thursday, July 7 by Zaanen, Bonsall, Amemiya, Wermer, 
Mikusi\'nski, Arens; on Friday, July 8 by Halperin, Gagliardo, Lax; on Monday, July 11 by Orlicz, Phillips, K\"othe, 
Mackey, Stampacchia, Nirenberg, Agmon; on Tuesday, July 12 by Nachbin, Helson, Kahane, Luxemburg and 
Fullerton.

Lectures from this conference appeared as proceedings [2], where we find the following 32 important contributions 
(with pages in [2]): 

{\small S. Agmon, {\it Remarks on self-adjoint and semi-bounded elliptic boundary value problems} (1-13), 
I. Amemiya, {\it On ordered topological linear spaces} (14-23), R. Arens, {\it The analytic-functional calculus in 
commutative Banach algebras} (24-28), N. Aronszajn, {\it Quadratic forms on vector spaces} (29-87), L. Bers, 
{\it Completeness theorems for Poincar\'e series in one variable} (88-100), F. F. Bonsall, {\it Semi-algebras of 
continuous functions} (101-114), J. Dieudonn\'e, {\it  Quasi-hermitian operators} (115-122), A. Dvoretzky, 
{\it Some results on convex bodies and Banach spaces} (123-160), L. Ehrenpreis, {\it A fundamental principle 
for systems of linear differential equations with constant coefficients, and some of its applications} (161-174), 
G. Fichera, {\it Spazi lineari di $k$-misure e di forme differenziali} (175-226), R. E. Fullerton, {\it Geometrical 
characterizations of certain function spaces} (227-236), E. Gagliardo, {\it A unified structure in various families 
of function spaces. Compactness and closure theorems} (237-241), I. Halperin, {\it Function spaces} (242-250), 
H. Helson and D. Lowdenslager, {\it Invariant subspaces} (251-262), E. Hille, {\it Linear differential equations 
in Banach algebras} (263-273), J.-P. Kahane, {\it Fonctions pseudo-p\'eriodiques dans $R^p$} (274-281), 
V. Klee, {\it Relative extreme points} (282-289), G. K\"othe, {\it Probleme der linearen Algebra in topologischen 
Vektorr\"aumen} (290-298), P. D. Lax, {\it Translation invariant spaces} (299-306), W. A. J. Luxemburg, 
{\it On closed linear subspaces and dense linear subspaces of locally convex topological linear spaces} 
(307-318), G. W. Mackey, {\it Induced representations and normal subgroups} (319-326), J. L. Massera, 
{\it Function spaces with translations and their application to linear differential equations} (327-334), J. Mikusinski, 
{\it Operations on distributions} (335-339), L. Nachbin, {\it Some problems in extending and lifting continuous 
linear transformations} (340-350), L. Nirenberg, {\it Inequalities in boundary value problems for elliptic differential 
equations} (351-356), W. Orlicz, {\it On spaces of $\phi $-integrable functions} (357-365), R. S. Phillips, 
{\it The extension of dual subspaces invariant under an algebra} (366-398), G. Stampacchia, {\it R\'egularisation 
des solutions de probl\`emes aux limites elliptiques \`a donn\'ees discontinues} (399-408), M. H. Stone, 
{\it Hilbert space methods in conformal mapping} (409-425), A. E. Taylor, {\it Spectral theory and Mittag-Leffler 
type expansions of the resolvent} (426-440), J. Wermer, {\it  Subalgebras of $C(X)$} (441-447), A. C. Zaanen, 
{\it Banach function spaces} (448-452).}
\vspace{3mm}

It should be mention that this photo was published in the book [1, p. 331] with the wrong information: {\it Robinson at 
the International Congress for Logic, Methodology, and Philosophy of Science, Jerusalem, 1964}, and without description 
of the names of persons on it. I have got this photo from the wife of W. Orlicz, after his death in 1990. He was my 
supervisor of Ph. D. defended in 1979 at the University of Pozna\'n in Poland.
\vspace{3mm}

{\bf Acknowledgments.} I wish to thank Dany Leviatan (Tel Aviv), Dan Amir (Tel Aviv), Jean-Pierre Kahane (Paris), Louis 
Nirenberg (New York), Edmond Ernest Granirer (Vancouver), Daniel Sternheimer (Dijon) and Shmuel Agmon (Jerusalem) 
for the help in recognition of persons on the photo and completing their dates of birth and death.

\vspace{3mm}

\centerline{ \bf References}
\vspace{2mm}

\noindent
[1] J. W. Dauben, {\it Abraham Robinson. The Creation of Nonstandard Analysis. A Personal 

and Mathematical Odyssey}, Princeton University Press, Princeton, NJ, 1995.

\noindent
[2] {\it Proceedings of the International Symposium on Linear Spaces} (Jerusalem, 1960), 

Jerusalem Academic Press, Pergamon Press, Oxford, Jerusalem 1961.

\vspace{7mm}

\noindent
{\small
Department of Engineering Sciences and Mathematics\\
Lule{\aa} University of Technology\\
971 87 Lule{\aa}, Sweden\\
e-mail: {\tt lech.maligranda@ltu.se}
}

\end{document}